\documentclass[12pt]{amsart}
\usepackage[cmtip,arrow]{xy}
\usepackage{
amsxtra,amsmath,amsthm,amssymb,amscd,ascmac,amsfonts,
multicol,tabularx,latexsym,mathrsfs}
\usepackage{amsfonts}
\theoremstyle{definition}

\newcommand{\Q}{\mathbb Q}
\newcommand{\Z}{\mathbb Z}

\newcommand{\Gal}{\mathrm{Gal}}
\newtheorem{thm}{Theorem}
\newtheorem{thmm}{Theorem}
\newtheorem{lem}{Lemma}
\newtheorem{prop}{Proposition}

\newtheorem{cor}{Corollary}

\newtheorem{ex}{Example}

\newcommand{\F}{\mathbb F}

\newcommand{\Cl}{\mathrm{Cl}}
\newdimen\minCDarrowwidth
\date{}
\begin{document}
\title[]
{On the $p$-adic limit of class numbers along a pro-$p$-extension}
\author{Manabu\ OZAKI}
\maketitle
\section{Introduction}
Let $p$ be a prime number and $K/k$ a $\Z_p$-extension over a number field $k$ of finite degree.

Sinnott \cite{S} found a mysterious phenomenon on the $p$-adic behavior
of the class numbers of the intermediate fields of $K/k$:
\begin{thmm}
Let $K/k$ be a $\Z_p$-extension over a number field $k$ of finite degree,
and we denote by $k_n$ the intermediate field of $K/k$
with $[k_n:k]=p^n\ (n\ge 0)$.
We denote by $h_n$ the class number of $k_n$.
Then the $p$-adic limit $\lim_{n\rightarrow\infty}h_n$ exists.
\end{thmm}
After Sinnott's work, Kisilevsky \cite{Ki} gave an alternative proof of the above theorem and generalized it to include function field case.
Recently, Ueki-Yoshizaki \cite{U-Y} rediscovered the above theorem and found similar phenomenon for $\Z_p$-covers of compact 3-manifolds and the torsion part of the first homology groups.  

In the present paper, we will first show that
a similar phenomenon occurs under much more general situation, namely,
for finitely generated pro-$p$-extensions (Theorem 1 in section 2).

The proof of Theorem 1 is based on the theory of representations
of finite groups.
However, in the case where $K/k$ is
the cyclotomic $\Z_p$-extension over
an abelian number field $k$, we will take an analytic approach to the theorem using $p$-adic $L$-functions and the analytic class number formula in section 3 (Theorem 7).

As the present paper was being written, we discovered the
earlier work of Han \cite{Han}, who establish
a similar result to Theorem 7 in section 3 for general CM-fields by using similar method to ours.

In section 4, based on the analytic method developed in section 3,
we will reveal enigmatic relationships between $p$-adic limits of various arithmetic invariants associated to $K/k$,
namely, the class numbers, the ratios of $p$-adic regulators and the square roots of discriminants, and the order of 
algebraic $K_2$-groups of the rings of integers.

\section{$p$-adic limits of class numbers in finitely generated pro-$p$-extensions}

Let $p$ be a prime number and 
$K/k$ a pro-$p$-extension over a number field
$k$ of finite degree whose Galois group is
a finitely generated pro-$p$-group.

We will show that the non $p$-parts of the class numbers of the intermediate fields of $K/k$ have the $p$-adic limit.
Specifically, we will show the following: 
\begin{thm}\label{thmclass}
Let $K/k$ be a pro-$p$-extension over a number field $k$ of finite degree whose Galois group is finitely generated.
We denote by $A(F)$ the non-$p$-part of the class group of $F$ for 
a number field $F$. 
Then, for any sequence
\[k_0\subseteq k_1\subseteq k_2\subseteq\cdots\subseteq k_i\subseteq k_{i+1}\subseteq\cdots
\]
of intermediate fields of $K/k$ so that $\bigcup_{i\ge 0} k_i=K$ and
$k_i/k$ is a finite normal extension,
there exists the $p$-adic limit $\lim_{n\rightarrow\infty}\#A(k_n)$.
Furthermore this limit does not depend on the
choice of the sequence $\{k_i\}$. 
\end{thm}
Here we note that Sinnott's theorem in Introduction follows from
the above theorem:

If $K/k$ is a $\Z_p$-extension, then the $p$-part of the class number
of the $n$-th layer converges to $0$ $p$-adically, or
stabilizes for all sufficiently large $n\ge 0$ by Iwasawa's class number formula.
Hence the class number of $k_n$ converges $p$-adically as $n\rightarrow\infty$ by Theorem \ref{thmclass}.

\

The above theorem follows from the following group theoretic general fact:
\begin{thm}\label{thmgrep}
Let $G$ be a finitely generated pro-$p$-group and
$A$ be a discrete torsion $G$-module whose $p$-primary part is trivial. We assume that $A^H$ is finite
for every open subgroup $H$ of $G$.
Then for any descending sequence
\[
H_0\supseteq H_1\supseteq\cdots\supseteq H_i\supseteq H_{i+1}\supseteq
\cdots\]
of open normal subgroups of $G$ with $\bigcap_{i\ge 0}H_i=1$, 
there exists the $p$-adic limit
$\lim_{n\rightarrow\infty}\# A^{H_n}$. 
Furthermore the above limit does not depend on the
choice of the descending sequence $\{H_i\}$.
\end{thm}
In fact, by applying Theorem 2 for $G=\Gal(K/k)$, $H_i=\Gal(K/k_i)$
and $A=A(K)$,
we get Theorem 1 because $A(K)^{\Gal(K/k_n)}\simeq A(k_n)$
holds:

For intermediate fields $F_1$ and $F_2$ of $K/k$ with
$k\subseteq F_1\subseteq F_2\subseteq K$,
let $j_{F_2/F_1}:A(F_1)\rightarrow A(F_2)$ be the natural lifting map
induced by the inclusion $F_1\subseteq F_2$.
Then $A(K)$ is the inductive limit $\varinjlim_{i\ge 0} A(k_i)$
of $A(k_i)$'s with respect to the maps $j_{k_i/k_j}$'s
($i\ge j\ge 0$).
Furthermore, if we denote by $N_{k_m/k_n}$ the norm map
from $A(K_m)$ to $A(k_n)$ ($m\ge n\ge 0$),
then we see that $N_{k_m/k_n}\circ j_{k_m/k_n}$ equals
the $[k_m:k_n]$-multiplication map on $A(k_n)$, which 
is an automorphism since the order of the finite abelian group $A(k_n)$
is prime to the power $[k_m:k_n]$ of $p$.
Then we find that $N_{k_m/k_n}$ is surjective and $j_{k_m/k_n}$ is injective.
Hence we have 
\begin{align*}
0=H^0(\Gal(k_m/k_n),A(k_m))&=A(k_m)^{\Gal(k_m/k_n)}/j_{k_m/k_n}\circ N_{k_m/k_n}(A(k_m))\\
&=A(k_m)^{\Gal(k_m/k_n)}/j_{k_m/k_n}(A(k_n)),
\end{align*}
which implies 
\[
A(k_m)^{\Gal(k_m/k_n)}=j_{k_m/k_n}(A(k_n))\ \ \ (m\ge n\ge 0),
\]
because 
$\#\Gal(k_m/k_n)$ is prime to $\#A(k_m)$.
Thus we conclude that 
\begin{align*}
A(K)^{\Gal(K/k_n)}
=\varinjlim_{m\ge n} A(k_m)^{\Gal(k_m/k_n)}=\varinjlim_{m\ge n} j_{k_m/k_n}A(k_n)\simeq A(k_n).
\end{align*}

\

In what follows we will give a proof of Theorem \ref{thmgrep}.
Let $G$ be a finitely generated pro-$p$-group.
Define $N_n$ to be the intersection of all the open subgroups $U$ of $G$ such that $[G:U]\le p^{n}$ for $n\ge 0$.
Then $N_n$ is a open normal subgroup of $G$ because $G$ is finitely generated. 
We will show the following:
\begin{prop}\label{cong}
Let $A$ be a discrete torsion $G$-module whose
$p$-part is trivial.
Assume that $A^H$ is finite for any open subgroup $H$ of $G$.
If $H$ and $H'$ are open normal subgroups of $G$ such that
$H, H'\subseteq N_n$,
then 
\[
\# A^H\equiv\# A^{H'}\pmod{p^n}
\]
holds.
\end{prop}
One can easily derive Theorem 2 from Proposition 1:

For any given $n\ge 0$, there exists a number $m_0$ such
that $H_i\subseteq N_n$ for all $i\ge m_0$,
because $\bigcap_{i\ge 0}H_i=1$ implies that $\{H_i\}$ is a fundamental
system of open neighborhoods of the identity in $G$.
Hence it follows from Proposition 1 that
$|\# A^{H_i}-\# A^{H_j}|_p\le p^{-n}$ holds for all $i,j\ge m_0$,
$|*|_p$ being the normalized $p$-adic absolute value on $\Z_p$.
This implies $\{\# A^{H_i}\}_{i\ge 0}$ is a $p$-adic Cauchy sequence in $\Z_p$.
Thus we have shown that the $p$-adic limit of $\{\# A^{H_i}\}_{i\ge 0}$ exists in $\Z_p$.

Furthermore, for any two descending
sequences $\{H_i\}$ and $\{H'_i\}$ of open normal subgroups with
$\bigcap_{i\ge 0}H_i=\bigcap_{i\ge 0}H'_i=1$ and $n\ge 0$,
if we choose $m_0$ so that $H_i,H'_i\subseteq N_n$ holds for every $i\ge m_0$,
then we have $|\# A^{H_i}-\# A^{H'_i}|\le p^{-n}$ for $i\ge m_0$.
Therefore the uniqueness of the limit also holds; 
$\lim_{i\rightarrow\infty}\# A^{H_i}
=\lim_{i\rightarrow\infty}\# A^{H'_i}$.

\

Now we will give a proof of Proposition 1.
It is enough to prove the proposition for
the case where $A$ is $l$-primary, 
$l\ne p$ being a prime number.

Furthermore we may assume that $H'\subseteq H$ holds:
Note that $H\cap H'\subseteq N_n$ is a open normal subgroup of $G$
for any open normal subgroups $H$ and $H'$ of $G$ with
$H, H'\subseteq N_n$. 
Then $\# A^{H\cap H'}\equiv\# A^{H}\pmod{p^n}$
and $\# A^{H\cap H'}\equiv\# A^{H'}\pmod{p^n}$
imply $\# A^{H}\equiv\# A^{H'}\pmod{p^n}$.

Put $\overline{G}:=G/H'\supseteq\overline{H}:=H/H'$, and
$B:=A^{H'}$.
Then $\overline{G}$ is a finite $p$-group, $B$ is a finite $l$-primary $\overline{G}$-module and
$B^{\overline{H}}\simeq A^{H}$.
Here we note that
$\overline{N_n}:=N_n/H'$ is the intersection of
all the subgroups of $\overline{G}$ whose indexes in $\overline{G}$
are less than or equal to $p^{n}$.

Therefore  it is enough to show that
$\# A\equiv \#A^H\pmod{p^n}$ if $H\subseteq N_n$
for the case where $G$ is a finite $p$-group and  
$A$ is a finite $l$-primary $G$-module.

For any $l$-primary finite $G$-module $X$ and $t\ge 0$, 
we define the $\F_l[G]$-module
$X(t)$ to be $l^tX/l^{t+1}X$.
Then we see that $X(t)^H\simeq X^H(t)$ holds for each $t\ge 0$
because $\# H$ is prime to $l$:

Since $\# l^{t+1}X$ is prime to $\# H$,
$H^1(H,l^{t+1}X)=0$ holds.
Hence we derive 
from the exact sequence
\[
0\longrightarrow l^{t+1}X\longrightarrow l^tX\longrightarrow X(t)\longrightarrow 0,
\]
the isomorphism 
$(l^tX)^H/(l^{t+1}X)^H\simeq X(t)^H$.
Furthermore, from the exact sequence
\[
0\longrightarrow X[l^t]\longrightarrow X\overset{l^t}{\longrightarrow}
l^tX\longrightarrow 0,
\]
we see $l^t(X^H)=(l^tX)^H$ by the similar reason.
Thus we have $X^H(t)\simeq X(t)^H$.

Hence $\# A(t)\equiv \# A(t)^H\pmod{p^n}$ for all $t\ge 0$
implies $\# A\equiv \# A^H\pmod{p^n}$.
Therefore we may assume that $A$ is a finite $\F_l[G]$-module.

Consequently, we have reduced a proof of Proposition 1
to that of the following:
\begin{prop}
Let $p$ and $l$ be distinct prime numbers.
Assume that $G$ is a finite $p$-group and $A$ a finite $\F_l[G]$-module.
Then, for any normal subgroup $H$ of $G$ and $n\ge 0$ with $H\subseteq N_n$,
$N_n$ being the intersection of all the subgroups $U$ of $G$ with
$[G:U]\le p^{n}$,
we have
\[
\# A\equiv \#A^H\pmod{p^n}.
\]
\end{prop}
\

\noindent
{\bf Proof of Proposition 2.}\ \ \ 
The $\F_l[G]$-module $A$ is isomorphic to the direct sum of
a finite number of irreducible $\F_l[G]$-modules;
\[
A\simeq\bigoplus_{\rho\in\mathrm{irr}_{\F_l}(G)} V_\rho^{e_\rho},
\]
holds for some non-negative integers $e_\rho$'s,
where $\mathrm{irr}_{\F_l}(G)$
and
$V_\rho$
denote the set of all the
irreducible representations of $G$ over $\F_l$ and 
the irreducible $\F_l[G]$-module associated to $\rho$, respectively.
Then we have 
\[
A^H\simeq\bigoplus_{\rho\in\mathrm{irr}_{\F_l}(G), 
\rho|_H=1}V_\rho^{e_\rho},
\]
since $\rho|_H\ne 1$ means $(V_\rho)^H=0$
by the normality of $H$ in $G$.
Thus we find that
\[
\#A/\#A^H=\prod_{\rho\in\mathrm{irr}_{\F_l}(G), 
\rho|_H\ne 1}(\# V_\rho)^{e_\rho}.
\]
Hence it is enough to show that
$\#V_\rho\equiv 1\pmod{p^n}$ for each 
$\rho\in\mathrm{irr}_{\F_l}(G)$ with $\rho|_H\ne 1$,
which follows from 
\begin{equation}\label{aim}
l^{\dim(\rho)}\equiv 1\pmod{p^n},
\end{equation}
where $\dim(\rho)$ stands for the dimension of $\rho$ over $\F_l$.

\

In what follows, we will show that \eqref{aim}
holds for $\rho\in\mathrm{irr}_{\F_l}(G)$
if $\rho|_H\ne 1$ and $H\subseteq N_n$.

We will employ the following theorem from theory of representations of finite $p$-groups to prove \eqref{aim}:
\begin{thm}\label{thmgprep}
Let $G$ be a finite $p$-group and 
$l$ a prime number different to $p$.
Let $\psi$ be 
an irreducible representation of $G$ over $\overline{\F_l}$,
an algebraic closure of $\F_l$.
Then there exist a subgroup $U$ of $G$ and a 1-dimensional representation
$\chi$ of $U$ over $\overline{\F_l}$ such that
$\psi=\mathrm{ind}_U^G(\chi)$ and either $\F_l(\psi)=\F_l(\chi)$
or $\F_l(\psi)\subseteq \F_l(\chi),\ [\F_l(\chi):\F_l(\psi)]=2$ holds
(the latter happens possibly for the case $p=2$).
Here we denote by $\F_l(\varphi)$ the field generated by the character values of
$\varphi$ over $\F_l$ for any representation $\varphi$ of $G$ over $\overline{\F_l}$.
\end{thm}
{\bf Proof.}\ \ \ The ``over $\Q$ version" of this theorem
is given in \cite[Theorem 1]{Ford87}.
We see that the over $\F_l$ version similarly follows because the characteristic of $\F_l$ is prime to $\# G$ by using 
a lifting of a representation from characteristic $l$ to characteristic 0.\hfill$\Box$

\

Let $\rho$ be any irreducible representation of $G$ over $\F_l$ with $\rho|_H\ne 1$,
and let $\psi$ be an irreducible factor of $\rho\otimes\overline{\F_l}$
over $\overline{\F_l}$.
Then we have
\begin{equation}\label{dec}
\rho\otimes\overline{\F_l}=\bigoplus_{\sigma\in\Gal(\F_l(\psi)/\F_l)
}{}^\sigma\psi,
\end{equation}
because, in general, the minimal splitting field of
a representation $\varphi$ of $G$ over a finite field of characteristic $l$
is exactly equal to $\F_l(\varphi)$.

On the other hand, it follows from Theorem \ref{thmgprep} that
there exist a subgroup $U$ of $G$ and 1-dimensional representation
$\chi$ of $U$ over $\overline{\F_l}$ such that $\psi=\mathrm{ind}_U^G(\chi)$.
Hence we see that
\begin{equation}\label{dim}
\dim(\rho)=\#\Gal(\F_l(\psi)/\F_l)\dim(\psi)=
[\F_l(\psi):\F_l][G:U]
\end{equation}
by \eqref{dec}.
\begin{lem}
(1)\ $\chi\ne 1$ holds.

(2)\ \ $\ker(\rho)=\bigcap_{g\in G} g\ker(\chi)g^{-1}$
\end{lem}

{\bf Proof.}\ \ \ 
(1)\ \ \ We see that $\mathrm{ind}_U^G(\overline{\F_l})\simeq\overline{\F_l}[G/U]$
is not an irreducible $\overline{\F_l}[G]$-module
if $U\ne G$, since $\overline{\F_l}[G/U]^G\simeq\overline{\F_l}$.
Because $\psi=\mathrm{ind}_U^G(\chi)$ is irreducible over
$\overline{\F_l}$, $\chi\ne 1$ if $U\ne G$.
If $U=G$, $\chi=1$ implies $\rho=1$, which contradicts to the assumption
$\rho|_H\ne 1$.
Thus we conclude that $\chi\ne 1$.

(2)\ \ \ It follows from \eqref{dec} that $\ker(\rho)=\ker(\psi)
=\ker(\mathrm{ind}_U^G(\chi))$.
Let $V$ be a representation space of $\chi$ and 
$G=\bigcup_{i=1}^d\sigma_iU$ a left coset decomposition of $G$
by $U$ ($d=[G:U]$).
Then a representation space of $\psi$
is given by $W:=V^d$ on which $\sigma\in G$ operates
\[
\sigma(x_i)_{1\le i\le d}=(\tau_{\pi^{-1}(i)}x_{\pi^{-1}(i)})_{1\le i\le d},\ \ 
((x_i)_{1\le i\le d}\in W)\]
where $\pi\in\frak{S}_d$ and $\tau_i\in U$ are given by
$\sigma\sigma_i=\sigma_{\pi(i)}\tau_i\ (1\le i\le d)$.
Hence we see that for $1\le i\le d$,
$\sigma\in G$ fixes $(0,\dots, 0,\overset{i}{x},0 \dots,0)\in W$
for all $x\in V$ if and only if
$\sigma\in\sigma_i\ker(\chi)\sigma_i^{-1}$.
Therefore $\sigma\in G$ fixes all the elements of $W$
if and only if $\sigma\in\bigcap_{i=1}^d\sigma_i\ker(\chi)\sigma_i^{-1}
=\bigcap_{g\in G}g\ker(\chi)g^{-1}$ (note that $\ker(\chi)$ is normal in $U$).
This means the assertion.\hfill$\Box$

\

Assume that $[G:U]=p^s$ and $\#\chi(U)=p^t\ \ (s,t\ge 0)$.
Then it follows from Lemma 1(1) that $t\ge 1$ and
we have $[G:\ker(\chi)]=[G:U][U:\ker(\chi)]=p^{s+t}$.
If $s+t\le n$, we have $[G:\ker(\chi)]\le p^n$, hence we find by using Lemma 1(2) that 
$\ker(\rho)=\bigcap_{g\in G} g\ker(\chi)g^{-1}\supseteq N_n\supseteq H$,
which contradicts to $\rho|_H\ne 1$.
Therefore $s+t\ge n+1$ must hold.

As well known, $f:=[\F_l(\chi):\F_l]=[\F_l(\mu_{p^t}):\F_l]$ is equal to
the order of $l\mod{p^t}$ in $(\Z/p^t\Z)^\times$.

In the case where $p\ne 2$, we find that
\[
l^{\dim(\rho)}=(l^f)^{p^s}\equiv 1\pmod{p^{t+s}}
\]
by using \eqref{dim} and $\F_l(\psi)=\F_l(\chi)$,
since $l^f\equiv 1 \pmod{p^t}$ and $t\ge 1$.
Because $s+t\ge n+1$, we have
$l^{\dim(\rho)}\equiv 1\pmod{p^{n+1}}$ in this case.

We will consider the case where $p=2$. 
In this case, $[\F_l(\psi):\F_l]=f$ or $\frac{f}{2}$ holds.
Then we find that
\[
l^{[\F_l(\psi):\F_l]}\equiv 1\pmod{2^{\max\{1,t-1\}}}.
\]
Therefore we have
\[
l^{\dim(\rho)}=l^{[\F_l(\psi):\F_l]2^s}\equiv 1\pmod{2^{t+s-1}}.
\]
Thus we conclude that $l^{\dim(\rho)}\equiv 1\pmod{2^n}$ holds
also in this case because $s+t\ge n+1$, 
which completes the proof of Proposition 2.
\hfill$\Box$

\

Thus we have proved Proposition \ref{cong}, which implies  Theorems \ref{thmclass} and \ref{thmgrep}.
\section{Analytic approach}
In this section, we will concentrate on the special and the most fundamental pro-$p$-extensions, namely,
the cyclotomic $\Z_p$-extensions
over abelian number fields.
In this case, we can approach to the $p$-adic limit of class numbers
via the analytic class number formula.

Let $k$ be an imaginary abelian number field, and $K/k$ the 
cyclotomic $\Z_p$-extension over $k$.
We denote by $k_n$ the $n$-th layer of $K/k$, and
let $h_n^-$ be the relative class number of the CM-field $k_n$.
The analytic class number formula expresses $h_n^-$
in terms of special values  of Dirichlet $L$-functions:
\begin{thm}[Analytic class number formula]\label{acnf}
\[
h_n^-=Q_nw_n\prod_{\chi\in\Gal(k_n/\Q)^\wedge_\mathrm{odd}}
\frac{1}{2}L(0,\chi),
\]
where $\Gal(k_n/\Q)^\wedge$ stands for the group of complex Dirichlet characters of $\Gal(k_n/\Q)$, 
$\Gal(k_n/\Q)^\wedge_\mathrm{odd}$ is defined to be
the set of all the odd characters in $\Gal(k_n/\Q)^\wedge$,
$Q_n$ is the unit index of the CM-field $k_n$, 
$w_n$ is the number of the roots of unity in $k_n$,
and $L(s,\chi)$ is the Dirichlet $L$-function associated to the 
Dirichlet character $\chi$.
\end{thm}

Here we regard Dirichlet characters as finite image characters of $\Gal(\Q^{\mathrm{ab}}/\Q)$ via the natural
isomorphism $\hat{\Z}^\times\simeq\Gal(\Q^{\mathrm{ab}}/\Q)$
given by the reciprocity map, where $\Q^{\mathrm{ab}}$ and $\hat{\Z}$
denote the maximal abelian extension field of $\Q$ and
the pro-finite completion of the ring $\Z$.

In what follows, we will give an analytic proof of the fact that the $p$-adic limit $\lim_{n\rightarrow\infty}h_n^-$ exists based on the analytic class number
formula and the Kubota-Leopoldt $p$-adic $L$-functions,
though this fact follows from Theorem 2 
in section 2:

We apply Theorem 2 to $G=\Gal(K/k)$, $H_i=\Gal(K/k_i)$, and
the inductive limit $A$ of 
\[
\Cl_n^-(p'):=\mathrm{ker}(N_{k_n/k_n^+}:\Cl_n(p')\longrightarrow \Cl_n^+(p')),
\]
where $\Cl_n(p')$ and $\Cl_n^+(p')$ are the non-$p$-parts of the class groups of $k_n$ and 
the maximal real subfield $k_n^+$ of $k_n$, respectively,
with respect to the natural lifting maps $\Cl_m(p')\longrightarrow\Cl_m(p')\ (m\ge n\ge 0)$.
Then we see that $A^{\Gal(K/k_n)}=\Cl_n(p')^-$. 
Therefore the $p$-adic limit
$\lim_{n\rightarrow\infty}h_n^-(p')$ of the non-$p$-parts
$h^-_n(p')$ of the relative class numbers $k_n$'s exists by Theorem 2.

On the other hand, Iwasawa's class number formulae for $K/k$
and $K^+/k^+$ imply that there exists the $p$-adic limit
$\lim_{n\rightarrow\infty}h_n^-(p)$ of the $p$-parts
$h^-_n(p)$ of the relative class numbers of $k_n$'s.
Therefore the $p$-adic limit
$\lim_{n\rightarrow\infty}h_n^-=
\lim_{n\rightarrow\infty}h_n^-(p')
\lim_{n\rightarrow\infty}h_n^-(p)$
exists.

\

We first recall the Kubota-Leopoldt $p$-adic $L$-function.
Let $\chi$ be a $p$-adic Dirichlet character, namely,
a group homomorphism from $\Gal(F/\Q)$ to $\overline{\Q_p}^\times$
for a certain abelian number field $F$.
Let 
\[\omega:\Gal(\Q(\mu_{2p})/\Q)\longrightarrow\Z_p^\times
\]
be the Teichm\"{u}ller character, namely, the character defined by
$\sigma(\zeta)=\zeta^{\omega(\sigma)}$
for $\sigma\in\Gal(\Q(\mu_{2p})/\Q)$
and $\zeta\in\mu_{p}$ (if $p\ne 2$) or $\zeta\in\mu_4$ (if $p=2$).
For any even $p$-adic Dirichlet character $\chi$,
there exists the $p$-adic meromorphic (analytic if $\chi\ne 1$)
function $L_p(s,\chi)$ 
defined on
$s\in\Z_p$ such that
\[
L_p(1-n,\chi)=(1-\chi\omega^{-n}(p)p^{n-1})L(1-n,\chi\omega^{-n})
\]
for all integers $n\ge 1$, and
$s=1$ is the unique pole of $L_p(s,1)$.
Here we fix 
embeddings $\overline{\Q}\hookrightarrow\overline{\Q_p}$ and
$\overline{\Q}\hookrightarrow\mathbb{C}$,
and then we regard $p$-adic Dirichlet characters as
complex Dirichlet characters and vise versa,
and view $L(1-n,\chi\omega^{-n})$ as a number in $\overline{\Q_p}$
since it is algebraic over $\Q$.

By using the analytic class number formula we can express
$h_n^-$ in terms of the values of $p$-adic $L$-functions
at $s=0$:
\begin{prop}\label{pancf}
Let $k$ be an imaginary abelian number field
whose conductor is not divisible by $p^2$ if $p\ne 2$ and by $8$ if $p=2$. Furthermore, we assume that $\chi(p)\ne 1$ for
any odd $\chi\in\Gal(k/\Q)^\wedge$. Then we have
\begin{equation}
h_n^-=Q_nw_n\prod_{\chi\in\Gal(k/\Q)^\wedge_\mathrm{odd}}
(1-\chi(p))^{-1}
\prod_{\substack{\chi\in\Gal(k/\Q)^\wedge_{\mathrm{odd}}\\
\psi\in\Gal(\Q_n/\Q)^\wedge}}
\frac{1}{2}L_p(0,\chi\psi\omega),
\end{equation}
$\Q_n$ being the $n$-th layer of the cyclotomic $\Z_p$-extension over
$\Q$.
\end{prop}
{\bf Proof.}\ \ \ We see that
\begin{equation}
L_p(0,\chi\psi\omega)=
\begin{cases}
(1-\chi(p))L(0,\chi)\ \ (\mbox{if $\psi=1$})\\
L(0,\chi\psi)\ \ (\mbox{if $\psi\ne 1$})
\end{cases}
\end{equation}
since $\chi\psi(p)=0$ if $\psi\ne 1$.
Hence the assertion of the proposition follows from Theorem \ref{acnf}.
\hfill$\Box$

\

Here we remark that if $\chi(p)=1$ for a certain odd $\chi\in\Gal(k/\Q)^\wedge$, then the Iwasawa $\lambda^-$-invariant
of $K/k$ is positive, which implies $h_n^-\rightarrow 0$
$p$-adically as $n\rightarrow \infty$.

\

Now we will show that the following theorem on the limit of the product
of $p$-adic $L$-functions along a $\Z_p$-extension,
which plays a crucial role in what follows:
\begin{thm}\label{limitL}
Let $\chi$ be an even Dirichlet character whose conductor is not divisible by $p^2$ if $p\ne 2$ and by $8$ if $p=2$.
We denote by $d_\chi$ the residue degree of 
$\Q_p(\chi)/\Q_p$, $\Q_p(\chi)$ being $\Q_p(\mathrm{im}(\chi))$.
We put
\[
\delta=
\begin{cases}
0\ (\mbox{if $p\ne 2$}),\\
1\ (\mbox{if $p=2$}).
\end{cases}
\]
(1)\ \ \ If $\chi\ne 1$, the $p$-adic limit
\[
\mathcal{L}(\chi):=
\lim_{n\rightarrow\infty}
\prod_{\psi\in\Gal(\Q_{nd_\chi}/\Q)^\wedge}2^{-\delta}L_p(s,\chi\psi)
\in\Z_p[\mathrm{im}(\chi)]
\]
exists and does not depend on $s\in\Z_p$.

Furthermore, we have
\begin{align*}
\prod_{\sigma\in\Gal(\Q_p(\chi)/\Q_p)}\mathcal{L}&({}^\sigma\!\chi)
\\=
&\lim_{n\rightarrow\infty}
\prod_{\psi\in\Gal(\Q_{n}/\Q)^\wedge}
\prod_{\sigma\in\Gal(\Q_p(\chi)/\Q_p)}
2^{-\delta}L_p(s,({}^\sigma\!\chi)\psi)
\in\Z_p.
\end{align*}

(2)\ \ \ There exists $c_p\in\Z_p^\times$ such that
\[
\lim_{n\rightarrow\infty}
p^{n+1+\delta}\prod_{\psi\in\Gal(\Q_n/\Q)^\wedge}2^{-\delta}L_p(s,\psi)
=\frac{c_p}{s-1}
\]
for $s\in\Z_p-\{1\}$.
\end{thm}
Before starting with the proof, 
we introduce Iwasawa's expression  of 
$p$-adic $L$-functions, which we will employ in the proof:
\begin{thm}
Let $\chi$ be a $p$-adic even Dirichlet character whose conductor
is not divisible by $p^2$ if $p\ne 2$ and by $8$ if $p=2$.
Put $\kappa=1+p$ if $p\ne 2$ or $\kappa=5$ if $p=2$.

(1)\ \ \ In the case where $\chi\ne 1$, there exists
the power series $f_\chi(T)\in 2\Z_p[\chi][[T]]$, 
$\Z_p[\chi]\subseteq\overline{\Q_p}$ being the ring generated by the values of $\chi$
over $\Z_p$, such that
\[
L_p(s,\chi\psi)
=f_\chi(\psi(\kappa)^{-1}\kappa^s-1)
\]
holds for $s\in\Z_p$ and any $\psi\in\Gal(\Q_n/\Q)^\wedge\ (n\ge 0)$.

(2)\ \ \ In the case where $\chi=1$, there exists
$g(T)\in2\Z_p[[T]]^\times$ such that
\[
L_p(s,\psi)=
\frac{g(\psi(\kappa)^{-1}\kappa^s-1)}{\psi(\kappa)^{-1}\kappa^s-\kappa}
\]
holds for any $\psi\in\Gal(\Q_n/\Q)^\wedge\ (n\ge 0)$.
\end{thm}

Next we introduce Coleman's norm operator:
Let $\mathcal{O}$ be the integer ring of a
finite extension $F$ of $\Q_p$, and let $\pi\in\mathcal{O}$ be
a prime element of $\mathcal{O}$.
Then there exists the unique map, so called Coleman's norm operator,
\[
\mathcal{N}:\mathcal{O}[[T]]\longrightarrow
\mathcal{O}[[T]]
\]
such that
\[
\mathcal{N}(f)((1+T)^{p}-1)
=\prod_{\zeta\in\mu_{p}}
f(\zeta(1+T)-1).
\]
See \cite[Lemma 13.39]{Wa} for a proof (That lemma treats the case where $\mathcal{O}=\Z_p$, however the proof there works for general cases).
Here we note that
\begin{equation}\label{N^n}
\mathcal{N}^n(f)((1+T)^{p^n}-1)
=\prod_{\zeta\in\mu_{p^n}}
f(\zeta(1+T)-1)
\end{equation}
holds for $n\ge 1$.

\begin{lem}\ \ \ 
For any $f(T)\in\mathcal{O}[[T]]^\times$,
$
\lim_{n\rightarrow\infty}\mathcal{N}^{dn}(f)\in\mathcal{O}[[T]]^\times
$
exists with respect to $\pi$-adic topology on $\mathcal{O}[[T]]$,
where $d$ stands for the residue degree of $F/\Q_p$.
\end{lem}
{\bf Proof.}\ \ \ 
The case where $\mathcal{O}=\Z_p$ is given by \cite[Cor.13.46]{Wa}.
We will follow the proof given in \cite{Wa} for general cases.

We first show that if $f(T)\in\mathcal{O}[[T]]^\times$ satisfies
$f((1+T)^p-1)\equiv 1\pmod{\pi^k}$, then 
$f(T)\equiv 1\pmod{\pi^k}$ holds:

Assume that $f(T)=1+\pi^m g(T)$ for $\pi\nmid g(T)\in\mathcal{O}[[T]]$
and $m\ge 0$.
Then we have $f((1+T)^p-1)=1+\pi^m g((1+T)^p-1)\equiv 1\pmod{\pi^k}$,
which implies 
$\pi^m g((1+T)^p-1)\equiv 0\pmod{\pi^k}$.
Since 
$g(T)\not\equiv 0\pmod{\pi}$, we see that
$g((1+T)^p-1)\equiv g(T^p)\not\equiv 0\pmod{\pi}$.
Hence we conclude that $m\ge k$ and $f(T)\equiv 1\pmod{\pi^k}$.

Noting that $f((((1+T)^p-1)+1)^{p^{n-1}}-1)=f((1+T)^{p^{n}}-1)$,
we have
\begin{equation}\label{pto1}
f((1+T)^{p^n}-1)\equiv 1\!\!\pmod{\pi^k}\Longrightarrow
f(T)\equiv 1\!\!\pmod{\pi^k}\ (k,n\ge 1).
\end{equation}

For $f(T)\in\mathcal{O}[[T]]]^\times$, we obtain
\begin{align*}
\mathcal{N}^d(f)((1+T)^{p^d}-1)&=
\prod_{\zeta\in\mu_{p^d}}f(\zeta(1+T)-1)
\equiv
f(T)^{p^d}\equiv f(T^{p^d})\\
&\equiv
f((1+T)^{p^d}-1)
\pmod{1-\zeta_{p^d}},
\end{align*}
$\zeta_p^d$ being a generator of $\mu_{p^d}$,
in $\mathcal{O}[\mu_d][[T]]$ since $\#\mathcal{O}/\pi=p^d$.
The above congruence yields
\[
\mathcal{N}^d(f)((1+T)^{p^d}-1)\equiv
f((1+T)^{p^d}-1)
\pmod{\pi}
\]
because $\mathcal{N}^d(f)((1+T)^{p^d}-1),\ f((1+T)^{p^d}-1)
\in\mathcal{O}[[T]]$
Hence we have 
\[
\frac{\mathcal{N}^d(f)((1+T)^{p^d}-1)}{f((1+T)^{p^d}-1)}\equiv 1\pmod{\pi},
\]
from which we find that 
\begin{equation}\label{=1}
\frac{\mathcal{N}^d(f)(T)}{f(T)}\equiv 1\pmod{\pi},
\end{equation}
by using \eqref{pto1}.

Assume that $f(T)\equiv 1\pmod{\pi^k}$.
Then $f(T)=1+\pi^kg(T)$ for a certain $g(T)\in\mathcal{O}[[T]]$, and we obtain
\begin{align*}
\mathcal{N}(f)((1+T)^p-1)&=\prod_{\zeta\in\mu_p}f(\zeta(1+T)-1)\\
&\equiv(1+\pi^kg(T))^p\equiv 1\pmod{\pi^k(\zeta_p-1)},
\end{align*}
in $\mathcal{O}[\mu_p][[T]]$, which yields
\[
\mathcal{N}(f)((1+T)^p-1)\equiv 1\pmod{\pi^{k+1}},
\]
and in turn implies $\mathcal{N}(f)(T)\equiv 1\pmod{\pi^{k+1}}$
by \eqref{pto1};
\begin{equation}\label{ktok+1}
f(T)\equiv 1\!\!\pmod{\pi^{k}}\Longrightarrow
\mathcal{N}(f)(T)\equiv 1\!\!\pmod{\pi^{k+1}}.
\end{equation}
Because the operator $\mathcal{N}$ is multiplicative, it follows from
\eqref{=1} and \eqref{ktok+1} that
if $i\ge j\ge 0$ and $i-j\equiv 0\pmod{d}$ then
\[
\mathcal{N}^i(f)(T)-\mathcal{N}^j(f)(T)\equiv 0\pmod{\pi^{j+1}}.
\]
Therefore the sequence $\{\mathcal{N}^{dn}(f)(T)\}_{n\ge 0}$ converges in $\mathcal{O}[[T]]$ for $f(T)\in\mathcal{O}[[T]]^\times$
since $\mathcal{O}[[T]]$ is complete with respect to $\pi$-adic topology.
\hfill$\Box$

\

We put
$\mathcal{N}^{d\infty}(f):=\lim_{n\rightarrow\infty}\mathcal{N}^{dn}(f)$
for $f(T)\in\mathcal{O}[[T]]^\times$.
\begin{lem}
For any $f(T)\in\mathcal{O}[[T]]^\times$, we have
\[
\lim_{n\rightarrow\infty}\prod_{\zeta\in\mu_{p^{dn}}}
f(\zeta\kappa^s-1)
=\mathcal{N}^{d\infty}(f)(0).
\]
\end{lem}
{\bf Proof.}\ \ \ 
It follows from \eqref{N^n} that
\begin{equation}
\mathcal{N}^{dn}(f)((\kappa^s)^{p^{dn}}-1)
=\prod_{\zeta\in\mu_{p^{dn}}}
f(\zeta\kappa^s-1)
\end{equation}

Because 
$\mathcal{N}^{dn}(f)(T)\rightarrow
\mathcal{N}^{d\infty}(f)(T)\ (n\rightarrow\infty)$ $\pi$-adically in $\mathcal{O}[[T]]$ by Lemma 2 and
$(\kappa^s)^{p^{dn}}-1\rightarrow 0\ (n\rightarrow\infty)$
$p$-adically in $\Z_p$,
we find that 
$\mathcal{N}^{dn}(f)((\kappa^s)^{p^{dn}}-1)\rightarrow
\mathcal{N}^{d\infty}(f)(0)$ as $n\rightarrow\infty$ for any $s\in\Z_p$.
Thus we obtain
\[
\lim_{n\rightarrow\infty}\prod_{\zeta\in\mu_{p^{dn}}}
f(\zeta\kappa^s-1)=\mathcal{N}^{d\infty}(f)(0).
\]
\hfill$\Box$

\

Now we will give a proof of Theorem 5:

{\bf Proof of Theorem 5.}
We first treat the case of $\chi\ne 1$.
Then it follows from Iwasawa's expression of the $p$-adic $L$-functions (Theorem 6) that
\[
\prod_{\psi\in\Gal(\Q_{nd_\chi}/\Q)^\wedge}2^{-\delta}L_p(s,\chi\psi)
=\prod_{\zeta\in\mu_{p^{nd_\chi}}}
2^{-\delta}f_\chi(\zeta\kappa^s-1).
\]
If $2^{-\delta}f_\chi(T)\not\in\Z_p[\chi][[T]]^\times$, we see that
\[
\lim_{n\rightarrow\infty}\prod_{\zeta\in\mu_{p^{nd_\chi}}}
2^{-\delta}f_\chi(\zeta\kappa^s-1)=0
\]
because, for $f(T)\in\mathcal{O}[[T]]$,
$v_p(f(\zeta\kappa^s-1))\ge\min\{v_p(f(0)),v_p(\zeta-1)=\frac{1}{\varphi(p^n)} \}$ if $\mu_{p^n}=\langle\zeta\rangle$, 
$v_p$ being the $p$-adic valuation on $\overline{\Q_p}$
normalized by $v_p(p)=1$.

If $2^{-\delta}f_\chi(T)\in\Z_p[\chi][[T]]^\times$, then we obtain
\[
\lim_{n\rightarrow\infty}\prod_{\zeta\in\mu_{p^{nd_\chi}}}
2^{-\delta}f_\chi(\zeta\kappa^s-1)=\mathcal{N}^{d_\chi\infty}
(2^{-\delta}f_\chi)(0)
\]
for any $s\in\Z_p$ by Lemma 3.

Furthermore, since 
\[
\prod_{\sigma\in\Gal(\Q_p(\chi)/\Q_p)}
2^{-\delta}f_{{}^\sigma\!\chi}(T)=
\prod_{\sigma\in\Gal(\Q_p(\chi)/\Q_p)}
2^{-\delta}\,{}^\sigma\!(f_\chi)(T)\in\Z_p[[T]],
\]
the latter part of (1) also follows
by Lemma 3.

\

Next we treat the case where $\chi=1$.
It follows from Theorem 6 that
\begin{equation}\label{L_p}
p^{n+1+\delta}\prod_{\psi\in\Gal(\Q_{n}/\Q)^\wedge}2^{-\delta}L_p(s,\psi)
=p^{n+1+\delta}\prod_{\zeta\in\mu_{p^n}}
\frac{2^{-\delta}g(\zeta\kappa^s-1)}{(\zeta\kappa^s-\kappa)}.
\end{equation}
Since $2^{-\delta}g(T)\in\Z_p[[T]]^\times$, 
we find that the $p$-adic limit
\begin{equation}\label{b_p}
b_p:=\lim_{n\rightarrow\infty}\prod_{\zeta\in\mu_{p^n}}
2^{-\delta}g(\zeta\kappa^s-1)\in\Z_p^\times
\end{equation}
exists, which is independent to $s$, by Lemma 3.

On the other hand, we have
\begin{equation}\label{prolim}
\begin{aligned}
\prod_{\zeta\in\mu_{p^n}}(\zeta\kappa^s-\kappa)
&=(-1)^\delta\kappa^{p^n}(\kappa^{p^n(s-1)}-1)\\
&=(-1)^\delta\kappa^{p^n}((1+p^{n+1+\delta}u_n)^{s-1}-1)\\
&\equiv(-1)^\delta\kappa^{p^n}(s-1)p^{n+1+\delta}u_n\pmod{p^{2n}}
\end{aligned}
\end{equation}for $n\ge 1$,
where $p\nmid u_n\in\Z$ is defined by
\[
\kappa^{p^n}=(1+p^{1+\delta})^{p^n}=1+u_np^{n+1+\delta}\ \ \ (n\ge 0).
\]
We see that the $p$-adic limit $u:=\lim u_n\in\Z_p^\times$ exists:
\[
\frac{u_{n+1}}{u_n}=
\frac{\kappa^{p^{n+1}}-1}{p(\kappa^{p^{n}}-1)}
=\frac{1+\kappa^{p^n}+\kappa^{2p^n}+\cdots+\kappa^{(p-1)p^n}
}{p}
\equiv 1\pmod{ p^{n+\delta}}
\]
since $\kappa^{p^n}\equiv 1\pmod{p^{n+1+\delta}}$.
This means $u_m-u_n\equiv 0\pmod{p^n}$ if $m\ge n\ge 0$.
Hence $\{u_n\}$ is a $p$-adic Cauchy sequence and converges.

Thus we conclude that
\begin{equation}\label{lim}
\lim_{n\rightarrow\infty}
\frac{1}{p^{n+1+\delta}}\prod_{\zeta\in\mu_{p^n}}(\zeta\kappa^s-\kappa)
=(-1)^\delta u(s-1)
\end{equation}
by using \eqref{prolim}.

Therefore we obtain
\[
\lim_{n\rightarrow\infty}
p^{n+1+\delta}\prod_{\psi\in\Gal(\Q_{n}/\Q)^\wedge}2^{-\delta}L_p(s,\psi)
=\frac{(-1)^\delta b_pu^{-1}}{s-1}=\frac{c_p}{s-1}
\]
if we put $c_p:=(-1)^\delta b_pu^{-1}\in\Z_p^\times$
by \eqref{L_p}, \eqref{b_p} and \eqref{lim}.
\hfill$\Box$

\

Now we will give an analytic proof of the fact that
the $p$-adic limit $\lim_{n\rightarrow\infty}h_n^{-}$ exists:

Recall that $k$ is an imaginary abelian number field whose conductor
is not divisible by $p^2$ if $p\ne 2$ and by 8 if $p=2$, and that $h_n^-$
is the relative class number of the $n$-th layer $k_n$ of the cyclotomic
$\Z_p$-extension $K/k$. 

Assume that every odd Dirichlet character $\chi\in\Gal(k/\Q)^\wedge$
satisfies $\chi(p)\ne 1$. If this does not hold, then we find $\lim_{n\rightarrow\infty}h_n^-=0$ as mentioned above (the paragraph after Proposition 3).
\begin{thm}
Under the above situation, the following holds:

\

(1)\ \ \ If $\mu_{2p}\not\subseteq k$, then 
\[
\lim_{n\rightarrow\infty}h_n^-
=b_0\prod_{\chi\in\Gal(k/\Q)^\wedge_\mathrm{odd}}
(1-\chi(p))^{-1}\prod_{\chi\in\Gal(k/\Q)^\wedge_\mathrm{odd}}
\mathcal{L}(\chi\omega),
\]
where 
\[
a_0:=w_0\lim_{n\rightarrow\infty}Q_n\in\Z,\ 
b_0=\begin{cases}
\omega(2)^{-\frac{[k:\Q]}{2}}a_0\ \ \mbox{(if $p\ne 2$)},\\
a_0\ \ \mbox{(if $p=2$)}.
\end{cases}
\]
(2)\ \ \ If $\mu_{2p}\subseteq k$, then
\[
\lim_{n\rightarrow\infty}h_n^-=
-b_1c_p\prod_{\chi\in\Gal(k/\Q)^\wedge_\mathrm{odd}}
(1-\chi(p))^{-1}
\prod_{1\ne\chi\in\Gal(k/\Q)^\wedge_\mathrm{even}}\mathcal{L}(\chi),
\]
where $\Gal(k/\Q)^\wedge_\mathrm{even}$ stands for the even characters in $\Gal(k/\Q)^\wedge$, and
\[
a_1:=\frac{w_0}{p^{1+\delta}}\lim_{n\rightarrow\infty}Q_n
\in\Z,\ \ 
b_1=\begin{cases}
\omega(2)^{-\frac{[k:\Q]}{2}}a_1\ \ \mbox{(if $p\ne 2$)},\\
a_1\ \ \mbox{(if $p=2$)}.
\end{cases}
\]
\end{thm}
{\bf Proof.}\ \ \ 
For each $\chi\in\Gal(k/\Q)^\wedge$ and $n\ge 0$, we put
\[
P_n(\chi):=\prod_{\psi\in\Gal(\Q_n/\Q)^\wedge}
\prod_{\sigma\in\Gal(\Q_p(\chi)/\Q_p)}2^{-\delta}
L_p(0,{}^\sigma\!\chi\omega\psi )
\]

Then we have
\begin{equation}\label{h_n^-}
h_n^-=2^{(\delta-1)\frac{[k_n:\Q]}{2}}Q_nw_n
\prod_{\chi\in\Gal(k/\Q)^\wedge,\mathrm{odd}}
(1-\chi(p))^{-1}
\prod_{i=1}^r P_n(\chi_i),
\end{equation}
where $\{\chi_1,\chi_2,\dots,\chi_r\}\subseteq\Gal(k/\Q)^\wedge$
is chosen so that
\[
\Gal(k/\Q)^\wedge_\mathrm{odd}=\coprod_{i=1}^r\{{}^\sigma\!\chi_i|
\sigma\in\Gal(\overline{\Q_p}/\Q_p)\},
\]
by using Proposition 3.

If $p\ne 2$, 
\begin{equation}\label{omega}
2^{(\delta-1)\frac{[k_n:\Q]}{2}}=(2^{p^n})^{-\frac{[k:\Q]}{2}}
\longrightarrow\omega(2)^{-\frac{[k:\Q]}{2}}
\end{equation}
$p$-adically as $n\rightarrow\infty$.

In the case where $\mu_{2p}\not\subseteq k$,
we see that $w_n=w_0$ for all $n\ge 0$.
Furthermore, if $Q_{n_0}=2$ holds for a certain $n_0\ge 0$, 
then $Q_n=2$ for all $n\ge n_0$:

$Q_{n}=2$ is equivalent to the surjectivity of the map
$\phi_n:E_n\longrightarrow\mu(k_{n}),\ \varepsilon\mapsto J(\varepsilon)/\varepsilon$,
where $E_n$, $J$, and $\mu(k_n)$ stand for
the unit group of $k_{n}$, the complex conjugation
of the CM-field $k_n$, and the group of the roots of unity
in $k_n$ (see the proof of \cite[Thm.4.12]{Wa}).
Since $\mu(k_n)=\mu(k_{n_0})$ and $E_{n_0}\subseteq E_n$, the assertion follows.

Hence $w_nQ_n=w_0Q_n$ stabilizes for sufficiently large $n\ge 0$,
especially the $p$-adic limit $a_0:=w_0\lim_{n\rightarrow\infty}Q_n$
exists in this case.
Therefore it follows from \eqref{h_n^-} and \eqref{omega} that
\[
\lim_{n\rightarrow\infty}h_n^-
=b_0\prod_{\chi\in\Gal(k/\Q)^\wedge_\mathrm{odd}}
(1-\chi(p))^{-1}\prod_{\chi\in\Gal(k/\Q)^\wedge_\mathrm{odd}}
\mathcal{L}(\chi\omega),
\]
where 
\[
b_0=\begin{cases}
\omega(2)^{-\frac{[k:\Q]}{2}}a_0\ \ \mbox{(if $p\ne 2$)},\\
a_0\ \ \mbox{(if $p=2$)},
\end{cases}
\]
since $P_n(\chi_i)\rightarrow
\prod_{\sigma\in\Gal(\Q_p(\chi_i)/\Q_p)}
\mathcal{L}({}^\sigma\!(\chi_i\omega))$ as $n\rightarrow\infty$
by Theorem 5.

In the case where $\mu_{2p}\subseteq k$, 
we see that $w_n=p^{n}w_0$ holds for all $n\ge 0$
(we note that $\mu_{2p^2}\not\subseteq k$ by the assumption on $k$).
Furthermore we see that if $Q_{n_0}=1$ holds for a certain $n_0\ge 0$,
then $Q_{n}=1$ holds for all $n\ge n_0$:

Because the norm map $N_{n+1,n}:\mu(k_{n+1})\longrightarrow\mu(k_n)$
is surjective in this case and 
\[
N_{n+1,n}(\phi_{n+1}(\varepsilon))=\phi_n(N_{n+1,n}(\varepsilon))
\ \ \ (\varepsilon\in E_{n+1}),
\]
if $\phi_{n+1}$ is surjective, then $\phi_n$ is also surjective, which means $Q_n=1\Longrightarrow Q_{n+1}=1$.

Hence $\frac{w_nQ_n}{p^{n+1+\delta}}=\frac{w_0Q_n}{p^{1+\delta}}\in\Z$
stabilizes for sufficiently large $n\ge 0$, especially the $p$-adic
limit $a_1:=\frac{w_0}{p^{1+\delta}}\lim_{n\rightarrow\infty}Q_n$ exists.
Then it follows from \eqref{h_n^-}, \eqref{omega} and Theorem 5 that
\begin{equation*}\label{h}
\lim_{n\rightarrow\infty}h_n^-=-
b_1c_p\prod_{\chi\in\Gal(k/\Q)^\wedge_\mathrm{odd}}
(1-\chi(p))^{-1}
\prod_{1\ne\chi\in\Gal(k/\Q)^\wedge_\mathrm{even}}\mathcal{L}(\chi),
\end{equation*}
where 
\[
b_1=\begin{cases}
\omega(2)^{-\frac{[k:\Q]}{2}}a_1\ \ \mbox{(if $p\ne 2$)},\\
a_1\ \ \mbox{(if $p=2$)},
\end{cases}
\]
since $p^{n+1+\delta}P_n(\omega^{-1})\rightarrow -c_p$ as $n\rightarrow\infty$
by Theorem 5
(Note that $\omega\in\Gal(k/\Q)^\wedge$ in this case).

Thus we have finished the analytic proof of the existence of the $p$-adic limit $\lim_{n\rightarrow\infty}h_n^{-}$.
\hfill$\Box$
\section{Relationships between the $p$-adic limits of various arithmetic
invariants along $\Z_p$-extensions}
In this section we will give certain relationships between the $p$-adic limits of various arithmetic invariants along $\Z_p$-extensions.

\

In what follows, we assume that
$k$ is an imaginary abelian number field with
$\mu_{2p}\subseteq k$ whose conductor is not divisible by
$p^2$ if $p\ne 2$ and by 8 if $p=2$.
Let $K/k$ be the cyclotomic $\Z_p$-extension with the $n$-th layer
$k_n$.
Also we put $\delta:=0$ if $p\ne 2$ or $\delta:=1$ if $p=2$ as in the preceding section.

The following theorem gives a relationship between the $p$-adic limits of the relative class numbers and the residues of $p$-adic zeta functions at $s=1$ along the $\Z_p$-extension $K/k$:
\begin{thm}
Let $h_n^-$ be the relative class number of $k_n$ and
$\rho_n$ the residue of $p$-adic zeta function $\zeta_{p,\,k_n^+}(s)$
of the maximal real subfield $k_n^+$ of $k_n$ at $s=1$.
Assume that $\chi(p)\ne 1$ for any $\chi\in\Gal(k/\Q)^\wedge_{\mathrm{odd}}$.
Then the $p$-adic limits
$h_\infty^-:=\lim_{n\rightarrow\infty}h_n^-$
and 
\[
\tilde{\rho}_\infty:=\lim_{n\rightarrow\infty}
p^{n+1+\delta-\frac{[k:\Q]}{2}2^n\delta}\rho_n\in\Z_p
\]
exist, and they satisfies
\[
h_\infty^-=
-b_1\tilde{\rho}_\infty\prod_{\chi\in\Gal(k/\Q)^\wedge_\mathrm{odd}}(1-\chi(p))^{-1},
\]
where $a_1:=\frac{w_0}{p^{1+\delta}}\lim_{n\rightarrow\infty}Q_n$ and
\[
b_1:=\begin{cases}
\omega(2)^{-\frac{[k:\Q]}{2}}a_1\ \ \mbox{(if $p\ne 2$)},\\
a_1\ \ \mbox{(if $p=2$)}.
\end{cases}
\]

\end{thm}
{\bf Proof.}\ \ \ 
By Theorem 7 (2),  we have
\[
h_\infty^-=-
b_1c_p\prod_{\chi\in\Gal(k/\Q)^\wedge_\mathrm{odd}}
(1-\chi(p))^{-1}
\prod_{1\ne\chi\in\Gal(k/\Q)^\wedge_\mathrm{even}}\mathcal{L}(\chi).
\]

On the other hand, because 
\begin{equation}\label{zeta}
\begin{aligned}
\zeta_{p,\,k_n^+}(s)&=\prod_{\chi\in\Gal(k_n/\Q)^\wedge_\mathrm{even}}
L_p(s,\chi)\\
&=\prod_{\psi\in\Gal(\Q_n/\Q)^\wedge}L_p(s,\psi)
\prod_{1\ne\chi\in\Gal(k/\Q)^\wedge_\mathrm{even}}
\prod_{\psi\in\Gal(\Q_n/\Q)^\wedge}L_p(s,\chi\psi),
\end{aligned}
\end{equation}
it follows from Theorem 5 that
\[
\tilde{\rho}_\infty=\lim_{n\rightarrow\infty}p^{n+1+\delta}2^{-\frac{[k:\Q]}{2}2^n\delta}\rho_n=c_p
\prod_{1\ne\chi\in\Gal(k/\Q)^\wedge_\mathrm{even}}
\mathcal{L}(\chi)\in\Z_p.
\]
Therefore we obtain
\[
h_\infty^-=
-b_1\tilde{\rho}_\infty\prod_{\chi\in\Gal(k/\Q)^\wedge_\mathrm{odd}}(1-\chi(p))^{-1}.
\]
\hfill$\Box$

\

Here we remark that if there exists  $\chi\in\Gal(k/\Q)^\wedge_\mathrm{odd}$
with $\chi(p)=1$, then
we have $h^-_\infty=\tilde{\rho}_\infty=0$,
because $\chi\omega\in\Gal(k/\Q)^\wedge_\mathrm{even}$
and $\mathcal{L}(\chi\omega)=0$ hold for such $\chi$.
\begin{ex}
We will treat the case where $p\ne 2$ and $k=\Q(\mu_p)$.
Then we have
\[
h_\infty^-=-(-1)^{\frac{p^2-1}{8}}2\tilde{\rho}_\infty
=-(-1)^{\frac{p^2-1}{8}}2\lim_{n\rightarrow\infty}p^{n+1}\rho_n,
\]
since $\omega(2)^{\frac{[k:\Q]}{2}}
=\omega(2)^\frac{p-1}{2}=(-1)^{\frac{p^2-1}{8}}$, $w_0=2p$
and $Q_n=1$ for all $n\ge 0$.

If $p=2$ and $k=\Q(\mu_4)$, then we have
\[
h_\infty^-=-\tilde{\rho}_\infty
=-\lim_{n\rightarrow\infty}
2^{n+2-2^n}\rho_n,
\]
since $w_0=4$
and $Q_n=1$ for all $n\ge 0$.

\end{ex}

\

The $p$-adic class number formula (see \cite[Thm.5.24]{Wa}, or 
\cite{Co}) says that
\[
\rho_n=\prod_{\frak{p}_n^+|p}\left(1-\frac{1}{N_n(\frak{p}_n^+)}
\right)
\frac{2^{p^n\frac{[k:\Q]}{2}}R_n^+h_n^+}{2\sqrt{D_n^+}},
\]
where $\frak{p}_n^+$ runs over all the primes of $k_n^+$ lying over $p$, $N_n(\frak{p}_n^+)$ denotes the norm of the ideal $\frak{p}_n^+$ of $k_n^+$, and $R_n^+$, $h_n^+$ and $D_n^+$ stand for the $p$-adic regulator,
the class number and the discriminant of $k_n^+$, respectively.
Here the signature of $\frac{R_n^+}{\sqrt{D_n^+}}$ is chosen so that
the above equality holds.

We see that
$\prod_{\frak{p}^+_n|p}\left(1-\frac{1}{N_n(\frak{p}^+_n)}\right)
=\prod_{\frak{p}^+_0|p}\left(1-\frac{1}{N_0(\frak{p}^+_0)}\right)$
holds for all $n\ge 0$ by the assumption on $k$.

Furthermore it follows from the theorem
in Introduction that the $p$-adic limit
$h_\infty^+:=\lim_{n\rightarrow\infty}h_n^+$ exists.
Hence we obtain
\[
\tilde\rho_\infty
=
\frac{\omega(2)^{\frac{[k:\Q]}{2}}h_\infty^+}{2}
\prod_{\frak{p}^+_0|p}\left(1-\frac{1}{N_0(\frak{p}^+_0)}\right)
\lim_{n\rightarrow\infty}
p^{n+1}\frac{R_n^+}{\sqrt{D_n^+}}
\]
if $p\ne 2$, and
\[
\tilde\rho_\infty
=
\frac{h_\infty^+}{2}\prod_{\frak{p}^+_0|p}\left(1-\frac{1}{N_0(\frak{p}^+_0)}\right)\lim_{n\rightarrow\infty}
2^{n+2}\frac{R_n^+}{\sqrt{D_n^+}}
\]
if $p=2$.

Therefore we have the following corollary to Theorem 8:
\begin{cor}Under the same assumptions of Theorem 8, we get
\begin{align*}
&h_\infty^-=\\
&-\frac{a_1h_\infty^+}{2}
\prod_{\chi\in\Gal(k/\Q)^\wedge_\mathrm{odd}}(1-\chi(p))^{-1}
\prod_{\frak{p}^+_0|p}\left(1-\frac{1}{N_0(\frak{p}^+_0)}\right)
\lim_{n\rightarrow\infty}p^{n+1+\delta}\frac{R_n^+}{\sqrt{D_n^+}}
\end{align*}
\end{cor}

\

\begin{ex}
Let $k=\Q(\mu_{2p})$.
Then we see that
\[
h_\infty^-=-\frac{h_\infty^+}{2}\left(1-\frac{1}{p}\right)
\lim_{n\rightarrow\infty}p^{n+1+\delta}\frac{R_n^+}{\sqrt{D_n^+}}
\]
holds.
\end{ex}

\

Finally, we will consider the $p$-adic limit of the orders of the algebraic $K_2$-groups of the integer rings along the $\Z_p$-extension $K^+/k^+$.

Let $\mathcal{O}_n^+$ be the integer ring of $k_n^+$ and
$K_2(\mathcal{O}_n^+)$ the algebraic $K_2$-group of
$\mathcal{O}_n^+$.
Then $K_2(\mathcal{O}_n^+)$ is a finite abelian group and
the order of it is given by:
\begin{thm}[The Birch-Tate conjecture (proved by \cite{Wi} combined with \cite{Ko})]
If $F$ is totally real finite abelian extension field of $\Q$, then we have
\[
\# K_2(\mathcal{O}_F)=w_2(F)|\zeta_F(-1)|,
\]
where $\zeta_F(s)$ is the Dedekind zeta function of $F$
and $w_2(F)$ stands for the number of the roots of unity in $F(\sqrt{F})$.
\end{thm}
The following theorem gives a relationship between the $p$-adic limits of $\# K_2(\mathcal{O}_n^+)$ and $h_n^-$:
\begin{thm}
The $p$-adic limit $\lim_{n\rightarrow\infty}2^{(-[k^+:\Q]2^n+1)\delta}\#K_2(\mathcal{O}_n^+)\in\Z_p$ exists
and we have
\begin{align*}
\lim_{n\rightarrow\infty}2^{(-[k^+:\Q]2^n+1)\delta}&\#K_2(\mathcal{O}_n^+)
\\
&=\frac{(-1)^{(1-\delta)[k^+:\Q]} v_0h_\infty^-}{2^{1-\delta}b_1(1-p)}
\prod_{\chi\in\Gal(k/\Q)^\wedge_\mathrm{odd}}
(1-\chi(p)),
\end{align*}
where 
$a_1=\frac{w_0}{p^{1+\delta}}\lim_{n\rightarrow\infty}Q_n$,
\[
b_1=\begin{cases}
\omega(2)^{-\frac{[k:\Q]}{2}}a_1\ \ \mbox{(if $p\ne 2$)},\\
a_1\ \ \mbox{(if $p=2$)},
\end{cases}
\]
and $v_0:=\frac{w_2(k_0^+)}{p^{1+\delta}}\in\Z$.
\end{thm}

\

{\bf Proof.}\ \ \ 
We first note that
\begin{align*}
\zeta_{p,k_n^+}(-1)&=
\prod_{\chi\in\Gal(k_n^+/\Q)^\wedge}L_p(-1,\chi)\\
&=(1-p)\prod_{\chi\in\Gal(k_n^+/\Q)^\wedge}L(-1,\chi\omega^{-2})
=(1-p)\zeta_{k_n^+}(-1)
\end{align*}
holds for $n\ge 0$ since $\omega^2\in\Gal(k_n^+/\Q)^\wedge$
and $\chi\omega^{-2}(p)=0$ if $\omega^2\ne\chi\in\Gal(k_n^+/\Q)^\wedge
$.

Furthermore, we find that
\[
w_2(k_n^+)=p^nw_2(k_0^+)=p^{n+1+\delta}v_0,
\]
where $v_0:=\frac{w_2(k_0^+)}{p^{1+\delta}}\in\Z$.

It follows from the functional equation of Dedekind zeta functions that
\[
|\zeta_F(-1)|=(-1)^{[F:\Q]}\zeta_F(-1)
\]
for any totally real number field $F$.

Therefore, by using \eqref{zeta}, Theorems 5 and 9, we get
\begin{align*}
\lim_{n\rightarrow\infty}&2^{(-[k^+:\Q]2^n+1)\delta}
\#K_2(\mathcal{O}_n^+)=
\lim_{n\rightarrow\infty}2^{(-[k^+:\Q]2^n+1)\delta}w_2(k_n^+)|\zeta_{k_n^+}(-1)|\\
&=
\frac{(-1)^{[k^+:\Q]} v_0}{1-p}\lim_{n\rightarrow\infty}
p^{n+1+\delta+(-[k^+:\Q]2^n+1)\delta}\zeta_{p,k_n^+}(-1)\\
&=
\frac{(-1)^{[k^+:\Q]} v_0}{1-p}
\frac{-c_p}{2^{1-\delta}}\prod_{1\ne\chi\in\Gal(k/\Q)^\wedge_\mathrm{even}}\mathcal{L}(\chi)\in\Z_p.
\end{align*}

Hence, comparing with Theorem 7 (2), we obtain
\begin{align*}
\lim_{n\rightarrow\infty}2^{(-[k^+:\Q]2^n+1)\delta}&\#K_2(\mathcal{O}_n^+)\\
&=\frac{(-1)^{[k^+:\Q]} v_0h_\infty^-}{2^{1-\delta}b_1(1-p)}
\prod_{\chi\in\Gal(k/\Q)^\wedge_\mathrm{odd}}
(1-\chi(p)).
\end{align*}
\hfill$\Box$

\

\begin{ex}
Let $k=\Q(\mu_{2p})$.
In the cese where $p\ne 2$, 
we see that $a_1=1,\ b_1=(-1)^{\frac{p^2-1}{8}}$, and that $v_0=12$ if $p\ne 3$
and $v_0=4$ if $p=3$.
Hence we have
\[
\lim_{n\rightarrow\infty}\#K_2(\mathcal{O}_n^+)
=
\begin{cases}
\frac{(-1)^{\frac{p-1}{2}+\frac{p^2-1}{8}}}
{1-p}6h_\infty^-\ \ 
\mbox{(if $p\ne 3$)},\\
-h_\infty^-\ \ \mbox{(if $p=3$)}.
\end{cases}
\]

In the case where $p=2$, we see that
$a_1=b_1=1$, and $v_0=3$.
Hence we have
\[
\lim_{n\rightarrow\infty}2^{-2^n+1}\#K_2(\mathcal{O}_n^+)=
3h_\infty^-.
\]
\end{ex}

\

\par\medskip\noindent
\noindent
Manabu Ozaki,\par\noindent
Department of Mathematics,\par\noindent
School of Fundamental Science and Engineering,\par\noindent
Waseda University,\par\noindent
Ohkubo 3-4-1, Shinjuku-ku, Tokyo, 169-8555, Japan\par\noindent
e-mail:\ \verb+ozaki@waseda.jp+

\end{document}